\documentclass[a4paper,twoside,10pt]{amsart}
\usepackage{amsmath,amsfonts,amssymb,amsthm}
\newtheorem{theorem}{Theorem}[section]
\newtheorem{cor}{Corollary}[section]
\newtheorem{cond}{Conditions}[section]
\usepackage{mathrsfs,enumerate}
\usepackage{color,graphicx}
\usepackage[a4paper]{geometry}
\author[G. Nemes]{Gerg\H{o} Nemes}
\address{School of Mathematics, The University of Edinburgh, James Clerk Maxwell Building, The King's Buildings, Peter Guthrie Tait Road, Edinburgh EH9 3FD, UK}
\email{Gergo.Nemes@ed.ac.uk}
\thanks{The author's research was supported by a research grant (GRANT11863412/70NANB15H221) from the National Institute of Standards and Technology.}


\begin{document}

\title[An extension of Laplace's method]{An extension of Laplace's method}

\begin{abstract} Asymptotic expansions are obtained for contour integrals of the form
\[
\int_a^b  \exp \left(  - zp(t) + z^{\nu /\mu } r(t) \right)q(t)dt,
\]
in which $z$ is a large real or complex parameter, $p(t)$, $q(t)$ and $r(t)$ are analytic functions of $t$, and the positive constants $\mu$ and $\nu$ are related to the local behaviour of the functions $p(t)$ and $r(t)$ near the endpoint $a$. Our main theorem includes as special cases several important asymptotic methods for integrals such as those of Laplace, Watson, Erd\'elyi and Olver. Asymptotic expansions similar to ours were derived earlier by Dingle using formal, non-rigorous methods. The results of the paper also serve to place Dingle's investigations on a rigorous mathematical foundation. The new results have potential applications in the asymptotic theory of special functions in transition regions, and we illustrate this by two examples.
\end{abstract}

\maketitle

\section{Introduction and main results}\label{intro}

Laplace's method is one of the best-known techniques for developing asymptotic approximations for integrals. The origins of the method date back to Laplace \cite{Laplace1812}, who studied the estimation of integrals of the form
\begin{equation}\label{eq32}
\int_a^b  e^{- zp(t)}q(t)dt.
\end{equation}
Here $(a, b)$ is a real (finite or infinite) interval, $z$ is a large positive parameter, and the functions $p(t)$ and $q(t)$ are continuous. Laplace made the observation that the major contribution to the integral \eqref{eq32} should come from the neighbourhood of the point where $p(t)$ attains its smallest value. Observe that by subdividing the range of integration at the minima and maxima of $p(t)$, and by reversing the sign of $t$ whenever necessary, we may assume, without loss of generality, that $p(t)$ has only one minimum in $\left[a, b\right]$ occurring at $t = a$. With certain assumptions on $p(t)$ and $q(t)$, Laplace's result is
\[
\int_a^b  e^{- zp(t)}q(t)dt \sim q(a)e^{ - z p(a)} \sqrt {\frac{\pi}{2 z p''(a)}}.
\]
The sign $\sim$ is used to mean that the quotient of the left-hand side by the right-hand side approaches $1$ as $z \to +\infty$.  This formula is now known as Laplace's approximation.

The modern version of the method was formulated by Erd\'elyi \cite[pp. 38--39]{Erdelyi1956}, who applied Watson's lemma to obtain a complete asymptotic expansion for the integral \eqref{eq32}.

The asymptotic theory of integrals of type \eqref{eq32} is also well established when $z$ is complex and $p(t)$, $q(t)$ are analytic functions in a domain of the complex plane containing the path of integration $\mathscr{P}$, say, joining $a$ to $b$. A well-known method for obtaining asymptotic expansions for such integrals is the method of steepest descents (for a detailed discussion of this method, see, e.g., \cite[Ch. 1]{Paris2011} or \cite[pp. 84--103]{Wong2001}). This method requires the deformation of $\mathscr{P}$ into a specific path that passes through one or more saddle points of $p(t)$ such that the function $\Im(p(t))$ is constant on it. (Recall that $t_0$ is a saddle point of $p(t)$ iff $p'(t_0)=0$.) This new path is called the path of steepest descent. However, in many specific cases, the construction of such a path can be extremely complicated. This problem may be bypassed using Perron's method which -- by requiring some extra assumptions -- avoids the computation of the path of steepest descent, and provides an explicit expression for the coefficients in the resulting asymptotic series \cite{Perron1917}\cite[Theorem 4, p. 105]{Wong2001}. A direct adaptation of Erd\'{e}lyi's theorem to complex integrals was formulated by Olver in the important paper \cite{Olver1970}, where explicit bounds are also provided for the error terms associated with the expansion (see also \cite[Theorem 6.1, p. 125]{Olver1974}).

In Chapter 9 of his book \cite{Olver1974}, Olver studied possible extensions of Laplace's method to integrals in which the large parameter $z$ enters in a more general way. He considered
\begin{equation}\label{eq38}
\int_a^b \exp \left(  - zp(t) + z^\sigma r(t) \right)q(t)dt,
\end{equation}
where $(a, b)$ is a real finite interval, $z$ is a large positive parameter, $\sigma<1$ is a constant, the function $p(t)$ is continuously differentiable in $\left(a,b\right]$, and the real or complex functions $q(t)$ and $r(t)$ are continuous in $\left(a,b\right]$. He assumed further that $p'(t)$ is positive in $\left(a,b\right]$,
\[
p(t) - p(a) = p_0(t - a)^\mu   + \mathcal{O}\left( (t - a)^{\mu_1} \right),\quad p'(t) = \mu p_0(t - a)^{\mu  - 1}  + \mathcal{O}\left( (t - a)^{\mu_1  - 1} \right) \qquad (t\to a+),
\]
and
\[
q(t) = q_0(t - a)^{\lambda  - 1}  + \mathcal{O}\left( (t - a)^{\lambda_1  - 1}\right),\quad r(t) = r_0(t - a)^\nu+\mathcal{O}\left( (t - a)^{\nu_1}  \right) \qquad (t \in \left(a,b\right]),
\]
where $p_0 > 0$, $q_0$ and $r_0$ are non-zero, $\mu_1>\mu  > 0$, $\lambda_1>\lambda > 0$ and $\nu_1> \nu  \ge 0$. Under these conditions, Olver showed that
\[
\int_a^b \exp \left(  - zp(t) + z^\sigma r(t) \right)q(t)dt  \sim \frac{q_0}{\mu }\Gamma \left( \frac{\lambda }{\mu } \right)\frac{e^{ - zp(a)}}{(p_0 z)^{\lambda /\mu } }
\]
as $z\to +\infty$, provided $\sigma < \min (1,\nu /\mu )$ \cite[Theorem 2.1, p. 326]{Olver1974}, and
\begin{equation}\label{eq37}
\int_a^b \exp \left(  - zp(t) + z^\sigma  r(t) \right)q(t)dt \sim \frac{q_0}{\mu }\operatorname{Fi}\left( \frac{\nu }{\mu },\frac{\lambda }{\mu };\frac{r_0}{p_0^{\nu /\mu }} \right)\frac{e^{ - zp(a)}}{(p_0 z)^{\lambda /\mu }}
\end{equation}
as $z\to +\infty$, when $\sigma = \nu /\mu$ \cite[Theorem 4.1, p. 332]{Olver1974}. In \eqref{eq37}, we use
\begin{equation}\label{eq35}
\operatorname{Fi}\left( \alpha ,\beta ;x \right) := \int_0^{+\infty}  \exp \left( - t + x t^\alpha \right)t^{\beta  - 1} dt \qquad (0\leq \alpha<1,\; \Re(\beta)>0)
\end{equation}
to denote the so-called Fax\'{e}n integral \cite[p. 332]{Olver1974}. With the parenthetic conditions it converges at both limits and defines an entire function of $x$. For basic properties, see Appendix \ref{appendixa}. Observe that the term $z^\sigma r(t)$ has no effect on the leading order asymptotic behaviour of the integral \eqref{eq38} when $\sigma < \min (1,\nu /\mu )$.

Integrals of type \eqref{eq38}, especially in the case that $\sigma = \nu /\mu$, play a crucial role in the asymptotic analysis of special functions in transition regions (see, e.g., \cite[Subsec. 5.2 and exer. 5.1, pp. 334--336]{Olver1974} or \cite{Olver1952}). Therefore, it is surprising that neither the complex analogue nor extension to a complete asymptotic expansion of the approximation \eqref{eq37} appears to have been recorded in the literature, apart from some formal results of Dingle (see below). The main purpose of the present paper is to fill these gaps by generalizing Olver's complex analogue of Erd\'elyi's method to integrals of the form
\begin{equation}\label{eq39}
I(z) := \int_a^b  \exp \left(  - zp(t) + z^{\nu /\mu } r(t) \right)q(t)dt,
\end{equation}
in which the path $\mathscr{P}$, say, is a contour in the complex plane joining $a$ to $b$, $z$ is a real or complex parameter, $\mu>\nu\geq 0$ are constants, and $p(t)$, $q(t)$ and $r(t)$ are analytic functions of $t$. The branch of $z^{\nu /\mu }$ has phase $(\nu /\mu )\arg z$. By imposing certain conditions on the functions $p(t)$, $q(t)$, $r(t)$ and the contour $\mathscr{P}$, we show that $I(z)$ possesses a complete asymptotic expansion, whose leading order term agrees with the right-hand side of \eqref{eq37}. We also provide explicit formulae for the asymptotic expansion coefficients similar to those given by Perron in the case of integrals of type \eqref{eq32}.

We begin by introducing the essential notation, which follows closely the notation of Olver \cite{Olver1970}. Let $t_1$ and $t_2$ be any two points of the contour $\mathscr{P}$. The part of $\mathscr{P}$ that lies between $t_1$ and $t_2$ will be denoted by $(t_1,t_2)_\mathscr{P}$ when $t_1$ and $t_2$ are both excluded, and by $\left[t_1,t_2\right]_\mathscr{P}$ when $t_1$ and $t_2$ are both included. Likewise for $\left(t_1,t_2\right]_\mathscr{P}$ and $\left[t_1,t_2\right)_\mathscr{P}$. We also denote
\[
\varpi = \text{angle of slope of} \; \mathscr{P} \; \text{at} \; a = \lim (\arg(t-a)) \qquad (t\to a \; \text{along} \; \mathscr{P}).
\]

Our assumptions are as follows:

\begin{cond}\label{cond1}
\begin{enumerate}[(i)]
	\item The functions $p(t)$, $q(t)$ and $r(t)$ are independent of $z$, and single-valued and analytic in a domain $\mathbf{T}$.
	\item The path of integration $\mathscr{P}$ is independent of $z$, $a$ is finite, $b$ is finite or infinite, and $(a,b)_\mathscr{P}$ lies within $\mathbf{T}$. (Hence, either $a$ or $b$ or both may be boundary points of $\mathbf{T}$.)
	\item In the neighbourhood of $a$, the functions $p(t)$, $q(t)$ and $r(t)$ can be expanded in convergent series of the form
\[
p(t) = p(a) + \sum\limits_{n = 0}^\infty p_n (t - a)^{n + \mu } ,\quad q(t) = \sum\limits_{n = 0}^\infty q_n (t - a)^{n + \lambda  - 1} ,\quad r(t) = \sum\limits_{n = 0}^\infty  r_n (t - a)^{n + \nu } ,
\]
where $p_0 \ne 0$, $\mu  > \nu  \ge 0$, and $\Re (\lambda ) > 0$. When $\mu$, $\lambda$ or $\nu$ is not an integer (which can occur only when $a$ is on the boundary of $\mathbf{T}$), the branches of $(t-a)^\mu$, $(t-a)^\lambda$ and $(t-a)^\nu$ are specified by the relations
\[
\left( t - a \right)^\mu  \sim \left| t - a \right|^\mu  e^{i\mu \varpi } ,\quad \left( t - a\right)^\lambda \sim \left| t - a \right|^\lambda  e^{i\lambda \varpi } ,\quad \left( t - a \right)^\nu  \sim \left| t - a \right|^\nu  e^{i\nu \varpi } ,
\]
as $t\to a$ along $\mathscr{P}$, and by continuity elsewhere on the contour $\mathscr{P}$.
\item $z$ ranges along a ray or over an annular sector given by $\theta_1 \leq \theta \leq \theta_2$ and $|z|\geq Z$, where $\theta = \arg z$, $\theta_2 - \theta_1<\pi$, and $Z>0$. Furthermore, $I\left(Z e^{i\theta}\right)$ converges at $b$ both absolutely and uniformly with respect to $\theta$.
\item $\left| \arg \left( e^{i\theta } p(t) - e^{i\theta } p(a) \right) \right| \le \frac{\pi }{2} - \delta$ with a fixed (small) positive $\delta$ when $t\in (a,b)_\mathscr{P}$, and $\Re \left( e^{i\theta } p(t) - e^{i\theta } p(a) \right)$ is bounded away from zero uniformly with respect to $\theta$ as $t\to b$ along the contour $\mathscr{P}$.
\item $r(t) = \mathcal{O}\left(\left| p(t) - p(a) \right|\right)$ as $t\to b$ along $\mathscr{P}$.
\end{enumerate}
\end{cond}

We are now in a position to formulate our main result.

\begin{theorem}\label{thm1} Subject to the forgoing assumptions,
\begin{equation}\label{eq34}
I(z) \sim e^{ - zp(a)} \sum\limits_{n = 0}^\infty  \frac{1}{z^{(n + \lambda )/\mu } }\sum\limits_{m = 0}^n f_{n,m} \operatorname{Fi} \left( \frac{\nu }{\mu },\frac{n + \lambda  + m\nu }{\mu };\frac{r_0 }{p_0^{\nu /\mu } } \right) 
\end{equation}
as $|z|\to +\infty$ uniformly in the sector $\theta_1 \leq \arg z \leq \theta_2$. The coefficients $f_{n,m}$ are expressible in terms of the $p_n$'s, $q_n$'s and $r_n$'s, for example,
\[
f_{0,0}  = \frac{q_0 }{\mu p_0^{\lambda /\mu } },\quad f_{1,0}  = \left( \frac{q_1 }{\mu } - \frac{(\lambda  + 1)p_1 q_0 }{\mu ^2 p_0 } \right)\frac{1}{p_0^{(1 + \lambda )/\mu } },\quad f_{1,1}  = \left( \frac{r_1 }{\mu } - \frac{\nu p_1 r_0}{\mu ^2 p_0 } \right)\frac{q_0}{p_0^{(1 + \lambda  + \nu )/\mu }}.
\]
Higher coefficients can be computed using formulae \eqref{eq20}, \eqref{eq21} and \eqref{eq22} below. The branch of $z^{(n + \lambda )/\mu}$ in \eqref{eq34} is $\exp\left((n+\lambda)(\log|z|+i\theta)/\mu\right)$.
\end{theorem}

\paragraph{\emph{Remarks}} 
\begin{enumerate}
  \item Note that, unlike in \eqref{eq37}, the upper limit of integration in \eqref{eq39} can be infinite.
	\item Provided consistency is maintained, neither $\varpi$ nor $\theta$ need be constrained to the principal range $\left(-\pi,\pi\right]$.
	\item Conditions (v) and (vi) may be replaced by (v') $\Re\left(e^{i\theta } p(t) - e^{i\theta } p(a)\right)$ is positive when $t\in (a,b)_\mathscr{P}$, and is bounded away from zero uniformly with respect to $\theta$ as $t\to b$ along $\mathscr{P}$, and (vi') $\Re (e^{i(\nu /\mu )\theta} r(t)) = \mathcal{O}\left( \Re\left(e^{i\theta } p(t) - e^{i\theta } p(a)\right) \right)$ uniformly with respect to $\theta$ as $t\to b$ along $\mathscr{P}$. In the special case that $r(t)$ is identically zero, Theorem \ref{thm1} with Conditions (i)--(iv) and (v') reduces to the result of Olver \cite{Olver1970}.
	\item If $\Re (e^{i(\nu /\mu )\theta} r(t))$ is non-positive on $\mathscr{P}$, Condition (vi) can be omitted.
\end{enumerate}

Expansions similar to \eqref{eq34} for integrals of type \eqref{eq39} when $\mu  = 2,3$, $\nu  = 1$ or $\mu  = 3$, $\nu  = 2$ are given in the book by Dingle \cite[Ch. X and XI]{Dingle1973}. Dingle's results, however, rely on interpretive, rather than rigorous, methods. The present results also serve to place Dingle's formal investigations on a rigorous mathematical foundation.

In the following corollary, we present an alternative form of the asymptotic expansion \eqref{eq34} in the special (and common) case that $\mu \geq 2$ is an integer and $\nu=1$. (In particular, $a$ is a saddle point of $p(t)$ in this case.) The main advantage of this form, over the form in \eqref{eq34}, is that it requires the computation of only $\mu$ different values of Fax\'{e}n's integral.

\begin{cor}\label{thm2} Assume that Conditions (i)--(vi) hold with $\mu \geq 2$ an integer and $\nu=1$. Then
\begin{align*}
I(z) & \sim e^{ - zp(a)} \sum\limits_{n = 0}^{\mu-1}  \frac{1}{z^{(\left\lceil n/2 \right\rceil+\lambda)/\mu }}\operatorname{Fi}^{(n)} \left( \frac{1}{\mu },\frac{\lambda }{\mu };\frac{r_0}{p_0^{1/\mu }} \right)\sum\limits_{m = 0}^\infty  \frac{\widetilde{f}_{n,m}}{z^{m/\mu }} \\ & = e^{ - zp(a)} \sum\limits_{n = 0}^{\mu-1}  \frac{1}{z^{(\left\lceil n/2 \right\rceil+\lambda)/\mu }}\operatorname{Fi} \left( \frac{1}{\mu },\frac{\lambda +n}{\mu };\frac{r_0}{p_0^{1 /\mu }} \right)\sum\limits_{m = 0}^\infty  \frac{\widetilde{f}_{n,m}}{z^{m/\mu }} .
\end{align*}
as $|z|\to +\infty$ uniformly in the sector $\theta_1 \leq \arg z \leq \theta_2$. Here the coefficients $\widetilde{f}_{n,m}$ can be determined by the procedure given in Section \ref{proof2}, and the branches of $z^{(\left\lceil n/2 \right\rceil+\lambda)/\mu }$ and $z^{m/\mu }$ are $\exp\left((\left\lceil n/2 \right\rceil+\lambda)(\log|z|+i\theta)/\mu\right)$ and $\exp\left(i m \theta/\mu\right)$, respectively.
\end{cor}

The remaining part of the paper is structured as follows. In Section \ref{proof}, we prove the asymptotic expansion stated in Theorem \ref{thm1}. Section \ref{proof2} discusses the proof of the alternative expansion given in Corollary \ref{thm2}. In Section \ref{app}, we give two illustrative examples to demonstrate the applicability of our results. The paper concludes with a short discussion in Section \ref{disc}.

\section{Proof of Theorem \ref{thm1}}\label{proof}

In this section, we prove the asymptotic expansion given in Theorem \ref{thm1}. Throughout the proof great care is required in prescribing the branches of the multi-valued functions that appear. In light of this, we introduce the following convention: the value of the angle $\varpi_0=\arg p_0$ is not necessarily the principal one, but rather is chosen so as to satisfy the inequality
\begin{equation}\label{eq23}
\left|\varpi_0 + \theta + \mu \varpi \right| \leq \frac{\pi}{2}.
\end{equation}
This branch of $\arg p_0$ is then to be used in constructing all of the fractional powers of $p_0$ which occur. For example, $p_0^{1/\mu}$ means $\exp\left((\log|p_0|+i\varpi_0)/\mu\right)$. Since
\[
e^{i\theta } p(t) - e^{i\theta } p(a) \sim e^{i\theta }p_0(t-a)^\mu
\]
as $t\to a$ along $\mathscr{P}$ (Condition (iii)), and $e^{i\theta } p(t) - e^{i\theta } p(a)$ has nonnegative real part (Condition (v)), it is always possible to choose $\varpi_0$ uniquely in this way. Furthermore, since $\theta$ is restricted to lie within an interval of length less than $\pi$, the value of $\varpi_0$ which satisfies \eqref{eq23} is independent of $\theta$.

We introduce the new variables $w$ and $v$ by
\begin{equation}\label{eq2}
w^\mu = v = p(t) - p(a).
\end{equation}
The branches of $\arg w$ and $\arg v$ are determined by the requirement that
\begin{equation}\label{eq24}
\mu \arg w, \; \arg v \to \varpi_0 +\mu \varpi \qquad (t\to a \; \text{along} \; \mathscr{P}),
\end{equation}
and by continuity elsewhere. Again, it is to be understood that these branches of $\arg w$ and $\arg v$ are used in constructing all fractional powers of $w$ and $v$ which occur. Since, by Condition (v), $v$ and $w$ are non-zero in $(a,b)_\mathscr{P}$, the branches are specified uniquely on $\mathscr{P}$; moreover, $\mu \arg w = \arg v$ at every point of $\mathscr{P}$. Condition (v) and \eqref{eq23}, \eqref{eq24} imply 
\begin{equation}\label{eq5}
\left|\theta + \arg v\right|<\frac{\pi}{2} \qquad (t \in (a,b)_\mathscr{P}).
\end{equation}
Therefore, $v$ is confined to a single Riemann sheet as $t$ traverses the contour $\mathscr{P}$.

For small $\left|t-a\right|$, Condition (iii) and the binomial theorem gives
\[
w= p_0^{1/\mu}(t-a)\left( 1+\frac{p_1}{\mu p_0}(t-a)+\cdots \right).
\]
Thus $w$ is a single-valued analytic function of $t$ in a neighbourhood of $a$, and $dw/dt$ is non-zero at $a$. Application of the inverse function theorem for analytic functions shows that for sufficiently small $\rho$, the disc $\left|t-a\right|<\rho$ is mapped conformally on a domain $\textbf{W}$ containing $w=0$ (see, e.g., \cite[Theorems 3.1 and 3.2, pp. 86--87]{Markushevich1965}). Furthermore, if $w\in \textbf{W}$, then $t-a$ can be expanded in a convergent series
\begin{equation}\label{eq13}
t - a = \sum\limits_{n = 1}^\infty a_n w^n  = \sum\limits_{n = 1}^\infty a_n v^{n/\mu } ,
\end{equation}
in which the $a_n$'s are expressible in terms of the $p_n$'s, for example, 
\[
a_1  = \frac{1}{p_0^{1/\mu } },\quad a_2  =  - \frac{p_1 }{\mu p_0^{1 + 2/\mu }},\quad a_3  = \frac{(\mu  + 3)p_1^2  - 2\mu p_0 p_2 }{2\mu ^2 p_0^{2 + 3/\mu }}.
\]

Let $k$ be a finite point of $\left(a,b\right]_\mathscr{P}$ chosen independently of $z$ and sufficiently close to $a$ to assure that the disc $\left|w\right|\leq \left|p(k)-p(a)\right|^{1/\mu}$ is contained in $\textbf{W}$. Then $\left[a,k\right]_\mathscr{P}$ may be deformed to make its $w$-map a straight line. Transformation to the variable $v$ yields
\begin{gather}\label{eq1}
\begin{split}
\int_a^k  \exp \left( - zp(t) + z^{\nu /\mu } r(t)  \right)q(t)dt &  = e^{ - zp(a)} \int_a^k  \exp \left(   - z\left(  p(t) - p(a)  \right) + z^{\nu /\mu } r(t)  \right)q(t)dt  
\\ & = e^{ - zp(a)} \int_0^\kappa  \exp \left(  - zv + \frac{r_0 }{p_0^{\nu /\mu } }(zv)^{\nu /\mu }  \right)f(v,z)dv,
\end{split}
\end{gather}
where $\kappa = p(k)-p(a)$ and
\begin{gather}\label{eq14}
\begin{split}
f(v,z) := & \exp \left( z^{\nu /\mu } \left( r(t) - \frac{r_0 }{p_0^{\nu /\mu } }v^{\nu /\mu } \right) \right)q(t)\frac{dt}{dv} \\ = & \exp \left( z^{\nu /\mu } \left( r(t) - \frac{r_0 }{p_0^{\nu /\mu } }(p(t) - p(a))^{\nu /\mu }  \right) \right)\frac{q(t)}{p'(t)}.
\end{split}
\end{gather}
The path of integration in the second line of \eqref{eq1} also is a straight line.

By combining the expansions of $p(t)$ and $q(t)$ (Condition (iii)) with the series \eqref{eq13}, we obtain the expansion
\begin{equation}\label{eq3}
\frac{q(t)}{p'(t)} = \sum\limits_{n = 0}^\infty b_n v^{(n + \lambda )/\mu  - 1},
\end{equation}
valid for sufficiently small values of $\left|v\right|$. The coefficients $b_n$ can be expressed in terms of the $p_n$'s and $q_n$'s, in particular,
\begin{gather*}
b_0  = \frac{q_0 }{\mu p_0^{\lambda /\mu }},\quad b_1  = \left( \frac{q_1}{\mu } - \frac{(\lambda  + 1)p_1 q_0 }{\mu ^2 p_0 } \right)\frac{1}{p_0^{(1 + \lambda )/\mu } },
\\ b_2  = \left( \frac{q_2}{\mu } - \frac{(\lambda  + 2)(p_1 q_1  + p_2 q_0 )}{\mu ^2 p_0 } + \frac{(\lambda  + \mu  + 2)(\lambda  + 2)p_1^2 q_0 }{2\mu ^3 p_0^2 } \right)\frac{1}{p_0^{(2 + \lambda )/\mu } }.
\end{gather*}
To obtain a representation for the general term, we use \eqref{eq2} and \eqref{eq3} together with Cauchy's formula:
\begin{equation}\label{eq4}
b_n  = \frac{1}{2\pi i}\oint_{(0 + )} \frac{q(t)}{p'(t)w^{\lambda  - \mu } } \frac{dw}{w^{n + 1}} = \frac{1}{\mu }\frac{1}{2\pi i}\oint_{(a + )} \frac{q(t)}{(p(t) - p(a))^{(n + \lambda )/\mu }}dt,
\end{equation}
with the contours of integration being positively-oriented small loops surrounding $0$ and $a$, respectively. Thus
\begin{equation}\label{eq20}
b_n  = \frac{1}{\mu }\mathop{\operatorname{Res}}\limits_{t = a} \left[ \frac{q(t)}{(p(t) - p(a))^{(n + \lambda )/\mu } } \right] \qquad (n=0,1,2,\ldots)
\end{equation}
(compare \cite[eq. 2.3.18]{DLMF}). Similarly, for small $\left|v\right|$, we have
\begin{equation}\label{eq18}
r(t) - \frac{ r_0 }{p_0^{\nu /\mu } }(p(t) - p(a))^{\nu /\mu }  = r(t) - \frac{ r_0 }{p_0^{\nu /\mu } }v^{\nu /\mu } = 
\sum\limits_{n = 1}^\infty c_n v^{(n+\nu)/\mu } ,
\end{equation}
where the coefficients $c_n$ are expressible in terms of the $p_n$'s and $r_n$'s, for example,
\[
c_1 = \left( r_1  - \frac{\nu p_1 r_0}{\mu p_0 } \right)\frac{1}{p_0^{(1 + \nu )/\mu }},\quad
c_2  = \left( r_2  - \frac{(\nu  + 1)p_1 r_1  + p_2 r_0 \nu }{\mu p_0 } + \frac{(\mu  + \nu  + 2)\nu p_1^2 r_0 }{2\mu ^2 p_0^2 } \right)\frac{1}{p_0^{(2 + \nu )/\mu } }.
\]
A calculation analogous to \eqref{eq4} shows that in general
\begin{equation}\label{eq21}
c_n  = \frac{1}{\mu }\mathop{\operatorname{Res}}\limits_{t = a} \left[ \frac{p'(t)r(t)}{(p(t) - p(a))^{(n + \nu )/\mu  + 1} } \right] \qquad (n=1,2,3,\ldots).
\end{equation}
Consequently, for small $\left|v\right|$, the exponential factor in \eqref{eq14} has the convergent expansion
\begin{equation}\label{eq6}
\exp \left( z^{\nu /\mu } \sum\limits_{n = 1}^\infty  c_n v^{(n + \nu )/\mu }  \right) = \sum\limits_{n = 0}^\infty  \left( \sum\limits_{m = 0}^n \frac{\textbf{\textsf{B}}_{n,m} (c_1 , \ldots ,c_{n - m + 1} )}{m!}(zv)^{m\nu /\mu }  \right)v^{n/\mu } ,
\end{equation}
where $\textbf{\textsf{B}}_{n,m}$ denotes the so-called partial ordinary Bell polynomials (see \cite[Appendix]{Nemes2013}).

Combination of \eqref{eq14}, \eqref{eq3}, \eqref{eq18} and \eqref{eq6} yields
\begin{equation}\label{eq29}
f(v,z) = \sum\limits_{n = 0}^\infty  \left( \sum\limits_{m = 0}^n f_{n,m} (zv)^{m\nu /\mu }  \right)v^{(n + \lambda )/\mu  - 1} ,
\end{equation}
with
\begin{equation}\label{eq22}
f_{n,m}  = \frac{1}{m!}\sum\limits_{j = m}^n b_{n - j} \textbf{\textsf{B}}_{j,m} (c_1 , \ldots ,c_{j - m + 1} ) .
\end{equation}
The first few of these coefficients are given by
\[
f_{0,0}  = b_0 ,\quad f_{1,0}  = b_1 ,\quad f_{1,1}  = b_0 c_1 ,\quad f_{2,0}  = b_2 ,\quad f_{2,1}  =  b_0 c_2 + b_1 c_1,\quad f_{2,2}  = \frac{1}{2}b_0 c_1^2 .
\]

For any non-negative integer $N$, define $f_N(v,z)$ by the relations
\begin{equation}\label{eq30}
f(v,z) = \sum\limits_{n = 0}^{N-1}  \left( \sum\limits_{m = 0}^n f_{n,m} (zv)^{m\nu /\mu }  \right)v^{(n + \lambda )/\mu  - 1} + v^{(N + \lambda )/\mu  - 1} f_N(v,z).
\end{equation}
Then, if $k$ is sufficiently close to $a$,
\begin{equation}\label{eq8}
f_N (v,z) = \mathcal{O} (1) \exp \left( C_{k}(zv)^{\nu /\mu } \right),
\end{equation}
where the implied constant is independent of $z$ and $v$, and $C_{k}$ is an assignable constant which may depend on $k$ (see Appendix \ref{appendixb}). We rearrange the integral in the second line of \eqref{eq1} in the form
\begin{gather}\label{eq9}
\begin{split}
& \int_0^\kappa  \exp \left(  - zv + \frac{r_0 }{p_0^{\nu /\mu }}(zv)^{\nu /\mu } \right)f(v,z)dv \\ & = \sum\limits_{n = 0}^{N - 1} \frac{1}{z^{(n + \lambda )/\mu } }\left( \sum\limits_{m = 0}^n f_{n,m} \int_0^{\infty e^{i\arg \kappa } } \exp \left( - zv + \frac{r_0 }{p_0^{\nu /\mu } }(zv)^{\nu /\mu }  \right)(zv)^{(n + \lambda  + m\nu )/\mu } \frac{dv}{v}  \right) \\ & \quad\; - \varepsilon _{N,2} (z) + \varepsilon _{N,1} (z),
\end{split}
\end{gather}
where $\arg \kappa$ is determined by
\begin{equation}\label{eq7}
\left|\theta + \arg \kappa \right| <\frac{\pi}{2}
\end{equation}
(cf. \eqref{eq5}), and
\begin{gather*}
\varepsilon _{N,1} (z) := \int_0^\kappa \exp \left( - zv + \frac{r_0 }{p_0^{\nu /\mu } }(zv)^{\nu /\mu } \right)v^{(N + \lambda )/\mu  - 1} f_N (v,z)dv,
\\
\varepsilon _{N,2} (z) := \sum\limits_{n = 0}^{N - 1} \frac{1}{z^{(n + \lambda )/\mu }}\left( \sum\limits_{m = 0}^n f_{n,m} \int_\kappa ^{\infty e^{i\arg \kappa } } \exp \left(  - zv + \frac{r_0}{p_0^{\nu /\mu } }(zv)^{\nu /\mu } \right)(zv)^{(n + \lambda  + m\nu )/\mu } \frac{dv}{v}  \right) .
\end{gather*}
All integration paths are now straight lines.

In consequence of \eqref{eq5} and \eqref{eq7} we derive
\begin{equation}\label{eq10}
\int_0^{\infty e^{i\arg \kappa } } \exp \left(  - zv + \frac{r_0}{p_0^{\nu /\mu }}(zv)^{\nu /\mu } \right)(zv)^{(n + \lambda  + m\nu )/\mu } \frac{dv}{v}  = \operatorname{Fi} \left( \frac{\nu }{\mu },\frac{n + \lambda  + m\nu }{\mu };\frac{r_0 }{p_0^{\nu /\mu } } \right).
\end{equation}

For $\varepsilon _{N,2} (z)$ the substitution $t=zv$ produces
\[
\varepsilon _{N,2} (z) = \sum\limits_{n = 0}^{N - 1} \frac{1}{z^{(n + \lambda )/\mu }}\left( \sum\limits_{m = 0}^n f_{n,m} \int_{z\kappa }^{\infty e^{i(\theta  + \arg \kappa )} } \exp \left(  - t + \frac{r_0 }{p_0^{\nu /\mu }}t^{\nu /\mu } \right)t^{(n + \lambda  + m\nu )/\mu  - 1} dt  \right) .
\]
If $0\leq \alpha<1$, and $\beta$, $x$ are arbitrary (fixed) complex numbers, and $\left|\zeta\right|$ is large, then
\begin{equation}\label{eq66}
\int_\zeta ^{\infty} \exp \left(  - t + xt^\alpha \right)t^{\beta  - 1} dt  = \mathcal{O}(1)\exp \left(  - \zeta  + x\zeta ^\alpha \right)\zeta ^{\beta  - 1} ,
\end{equation}
uniformly in the region $\left|\arg \zeta\right| \leq \frac{\pi}{2}$ (see Appendix \ref{appendixc}). The integration path in \eqref{eq66} is a straight line starting at $\zeta$ and making a positive or negative acute angle with the real axis. Hence
\begin{equation}\label{eq11}
\varepsilon _{N,2} (z) = \mathcal{O}(1)\exp \left(  - z\kappa  + \frac{r_0}{p_0^{\nu /\mu }}(z\kappa )^{\nu /\mu } \right) z^{(N-1)\nu /\mu  - 1} ,
\end{equation}
uniformly with respect to $\theta$ as $\left|z\right| \to +\infty$.

For $\varepsilon _{N,1} (z)$ the substitution $v = \kappa \tau$ and the estimate \eqref{eq8} give
\[
\varepsilon _{N,1} (z) = \int_0^1 \exp \left(  - z\kappa \tau  + \left( {\frac{r_0}{p_0^{\nu /\mu } } + C_{k}} \right)(z\kappa \tau )^{\nu /\mu } \right)\tau ^{(N + \lambda )/\mu  - 1} \mathcal{O}(1)d\tau .
\]
In consequence of Condition (v) and the fact that $\theta$ is restricted to a closed interval, we have
\[
\Re (z\kappa ) = \left| z \right|\Re \left( {e^{i\theta } p(k) - e^{i\theta } p(a)} \right) \ge \left| z \right|\eta _k ,
\]
where $\eta_k$ is independent of $z$ and positive. Hence, we may assert that
\begin{align*}
\left| \varepsilon _{N,1} (z) \right| & \le \int_0^1 \exp \left(  - \left| z \right|\eta _k \tau  + \left| \frac{r_0 }{p_0^{\nu /\mu }} + C_k  \right|\left(\left| z \right| \left| \kappa  \right| \tau\right) ^{\nu /\mu }  \right)\tau ^{\left( N + \Re (\lambda ) \right)/\mu  - 1} \mathcal{O}(1)d\tau 
\\ & \le \frac{1}{\left| z \right|^{\left( N + \Re (\lambda ) \right)/\mu }}\int_0^{ + \infty } \exp \left(  - \eta _k t + \left| \frac{r_0}{p_0^{\nu /\mu }} + C_k  \right|\left| \kappa  \right|^{\nu /\mu } t^{\nu /\mu } \right)t^{\left( N + \Re (\lambda ) \right)/\mu  - 1} \mathcal{O}(1)dt .
\end{align*}
In arriving at the second inequality, we have made a change of integration variable from $\tau$ to $t$ via $t=\left| z \right| \tau$ and extended the upper limit of integration to $+\infty$. Since $0\leq \nu/\mu <1$, the integral in the second line is seen to be convergent, and thus
\begin{equation}\label{eq12}
\varepsilon _{N,1} (z) = \mathcal{O}\left( \left| z \right|^{ - (N + \Re (\lambda ))/\mu } \right) = \mathcal{O}\left( \left| z \right|^{ - (N + \lambda )/\mu } \right),
\end{equation}
uniformly with respect to $\theta$.

Combination of \eqref{eq1}, \eqref{eq9}, \eqref{eq10}, \eqref{eq11} and \eqref{eq12} yields
\begin{gather}\label{eq19}
\begin{split}
& \int_a^k  \exp \left( - zp(t) + z^{\nu /\mu } r(t)  \right)q(t)dt \\ & = e^{ - zp(a)} \left(\sum\limits_{n = 0}^{N-1}  \frac{1}{z^{(n + \lambda )/\mu } }\sum\limits_{m = 0}^n f_{n,m} \operatorname{Fi} \left( \frac{\nu }{\mu },\frac{n + \lambda  + m\nu }{\mu };\frac{r_0 }{p_0^{\nu /\mu } } \right) + \mathcal{O}\left( \frac{1}{\left|z\right|^{(N + \lambda )/\mu }} \right) \right),
\end{split}
\end{gather}
uniformly with respect to $\theta$ as $\left|z\right| \to +\infty$.

It remains to consider the tail of the original integral, that is, the contribution from $(k,b)_\mathscr{P}$. We choose a number $Z$ satisfying Condition (iv) and write 
\begin{align*}
& \int_k^b\exp \left(  - zp(t) + z^{\nu /\mu } r(t) \right)q(t)dt = e^{ - zp(a)} \int_k^b \exp \left(  - z\left( p(t) - p(a) \right) + z^{\nu /\mu } r(t) \right)q(t)dt
\\ & = e^{ - zp(a)} e^{Ze^{i\theta } p(a)} \int_k^b \exp \left(  - \left( \left| z \right| - Z \right)\left( e^{i\theta } p(t) - e^{i\theta } p(a) \right) + \left( \left| z \right|^{\nu /\mu }  - Z^{\nu /\mu }  \right)e^{i(\nu /\mu )\theta } r(t) \right) \\ & \qquad\qquad\qquad\qquad \times \exp \left(  - Ze^{i\theta } p(t) + \left( Ze^{i\theta }  \right)^{\nu /\mu } r(t) \right)q(t)dt .
\end{align*}
From Condition (v) it follows that
\[
\Re \left(e^{i\theta } p(t) - e^{i\theta } p(a)\right) \geq \eta >0, \qquad (t \in \left[k,b\right)_\mathscr{P}),
\]
where $\eta$ is independent of $\theta$. From Conditions (i), (ii) and (vi), we can assert that $\left| r(t) \right| \le K_1\left| p(t) - p(a) \right| +K_2$ with some appropriate positive constants $K_1$ and $K_2$ provided $t\in \left[k,b\right)_\mathscr{P}$. Accordingly,
\begin{align*}
& -\left(  \left| z \right| -Z\right)\Re \left(e^{i\theta } p(t) - e^{i\theta } p(a)\right) + \left( \left| z \right|^{\nu /\mu }  - Z^{\nu /\mu }  \right)\Re \left(e^{i(\nu /\mu )\theta} r(t)\right)
\\ & \le -\left(  \left| z \right|-Z\right)\Re \left(e^{i\theta } p(t) - e^{i\theta } p(a)\right) + \left( \left| z \right|^{\nu /\mu }  - Z^{\nu /\mu }  \right)\left| r(t) \right|
\\ & \le -\left( \left| z \right| -Z\right)\Re \left(e^{i\theta } p(t) - e^{i\theta } p(a)\right) + \left( \left| z \right|^{\nu /\mu }  - Z^{\nu /\mu }  \right)\left( K_1 \left| p(t) - p(a) \right| + K_2  \right)
\\ & \le -\left( \left| z \right| -Z\right)\Re \left(e^{i\theta } p(t) - e^{i\theta } p(a)\right) + \left( \left| z \right|^{\nu /\mu }  - Z^{\nu /\mu }  \right)\left( \frac{K_1}{\cos \delta }\Re \left(e^{i\theta } p(t) - e^{i\theta } p(a)\right) + K_2  \right)
\\ & = - \left( \frac{1}{2}\left( \left| z \right| - Z \right) - \left(\left| z \right|^{\nu /\mu }  - Z^{\nu /\mu } \right)\frac{K_1}{\cos \delta} \right)\Re \left(e^{i\theta } p(t) - e^{i\theta } p(a)\right) \\ & \quad \, - \frac{1}{2}\left( \left| z \right| - Z \right)\Re \left(e^{i\theta } p(t) - e^{i\theta } p(a)\right) + \left( \left| z \right|^{\nu /\mu }  - Z^{\nu /\mu } \right)K_2
\\ & \le - \left( \frac{1}{2}\left( \left| z \right| - Z \right) - \left( \left| z \right|^{\nu /\mu }  - Z^{\nu /\mu }  \right)\frac{K_1}{\cos \delta } \right)\eta  - \frac{1}{2}\left( \left| z \right| - Z\right)\eta  + \left( {\left| z \right|^{\nu /\mu }  - Z^{\nu /\mu } } \right)K_2  \le  - \varepsilon \left| z \right|,
\end{align*}
for some $\varepsilon>0$ and sufficiently large $|z|$. Consequently,
\begin{align*}
\left| \int_k^b \exp \left(  - zp(t) + z^{\nu /\mu } r(t) \right)q(t)dt \right| \le \; & \left| e^{ - zp(a)} \right|\left| e^{Ze^{i\theta } p(a)}  \right|e^{ - \varepsilon \left| z \right|} \\ & \times \int_k^b \left| \exp \left(  - Ze^{i\theta } p(t) + \left( Ze^{i\theta }  \right)^{\nu /\mu } r(t) \right)q(t)dt \right| .
\end{align*}
Condition (iv) shows that this upper bound is $e^{ - zp(a)} \mathcal{O}\left( e^{ - \varepsilon \left| z \right|} \right)$, uniformly with respect to $\theta$. Therefore, the asymptotic expansion \eqref{eq19} is unaffected by addition of the contribution from $(k,b)_\mathscr{P}$.

\section{Proof of Corollary \ref{thm2}}\label{proof2}

In this section, we prove the alternative expansion given in Corollary \ref{thm2}. Throughout we assume that $\mu \geq 2$ is an integer. It is readily seen from \eqref{eq35} that for any non-negative integer $n$,
\[
\operatorname{Fi}^{(n)} \left( \alpha ,\beta ;x \right) := \frac{d^n}{dx^n}\operatorname{Fi}\left( \alpha ,\beta ;x \right) = \operatorname{Fi}\left( \alpha ,\beta  + n\alpha ;x \right) .
\]
Thus, when $\nu=1$, the asymptotic expansion \eqref{eq34} may be written
\begin{equation}\label{eq36}
I(z) \sim e^{ - zp(a)} \sum \limits_{n = 0}^\infty  \frac{1}{z^{(n + \lambda )/\mu } }\sum\limits_{m = 0}^n f_{n,m} \operatorname{Fi}^{(n + m)} \left( \frac{1}{\mu },\frac{\lambda }{\mu };\frac{r_0}{p_0^{1 /\mu } } \right) .
\end{equation}
According to the results of the papers \cite{Kaminski1997,Paris1985}, the function
\begin{equation}\label{eq46}
w(x):=\operatorname{Fi}\left( \frac{1}{\mu },\frac{\lambda }{\mu };x \right)
\end{equation}
satisfies the $\mu$th order linear ordinary differential equation
\begin{equation}\label{eq44}
\mu w^{(\mu)} (x) = xw'(x) + \lambda w(x).
\end{equation}
Consequently, we may write
\begin{equation}\label{eq45}
w^{(j)} (x) = \sum\limits_{n = 0}^{\mu  - 1} g_{j,n} (x)w^{(n)} (x),
\end{equation}
where the coefficients $g_{j,n} (x)$ are polynomials in $x$. It is clear that when $j < \mu$,
\[
\left.\begin{aligned}
&  g_{j,j} (x) = 1,\\
&  g_{j,n} (x) = 0,\quad n \neq j.
\end{aligned}\right\}
\]
When $j \ge \mu$, we differentiate both sides of the equality \eqref{eq45} and employ \eqref{eq44}. This yields the recurrence scheme
\[
\left.\begin{aligned}
  &  g_{\mu ,0} (x) = \lambda/\mu,\; g_{\mu ,1} (x) = x/\mu,\; g_{\mu ,2} (x) =  \cdots  = g_{\mu ,\mu  - 1} (x) = 0, \\
  &  g_{j + 1,0} (x) = g'_{j,0} (x) + (\lambda/\mu)g_{j,\mu  - 1} (x),\\
  &  g_{j + 1,1} (x) = g_{j,0} (x) + g'_{j,1} (x) + (x/\mu)g_{j,\mu  - 1} (x),\\
  &  g_{j + 1,n} (x)  = g_{j,n - 1} (x) + g'_{j,n} (x),\quad 2 \leq n \leq \mu-1.
\end{aligned}\right\}
\]
Rearranging the terms in \eqref{eq36} and employing \eqref{eq45} (together with \eqref{eq46}), we deduce
\begin{align*}
I(z) & \sim e^{ - zp(a)} \sum\limits_{j = 0}^\infty \operatorname{Fi}^{(j)} \left( \frac{1}{\mu },\frac{\lambda }{\mu };\frac{r_0}{p_0^{1 /\mu } } \right)\sum\limits_{m = \left\lceil {j/2} \right\rceil }^j \frac{f_{m,j - m} }{z^{(m + \lambda )/\mu } } 
\\ &= e^{ - zp(a)} \sum\limits_{n = 0}^{\mu  - 1} \operatorname{Fi}^{(n)} \left( \frac{1}{\mu },\frac{\lambda }{\mu };\frac{r_0}{p_0^{1 /\mu } } \right)\sum\limits_{m = 0}^\infty  \frac{1}{z^{(m + \lambda )/\mu }}\sum\limits_{j = m}^{2m} g_{j,n} \left( \frac{r_0}{p_0^{1 /\mu } } \right)f_{m,j - m} 
\\ &= e^{ - zp(a)} \sum\limits_{n = 0}^{\mu  - 1} \frac{1}{z^{(\left\lceil n/2 \right\rceil + \lambda )/\mu } }\operatorname{Fi}^{(n)} \left( \frac{1}{\mu },\frac{\lambda }{\mu };\frac{r_0}{p_0^{1 /\mu } } \right) \\ & \qquad\qquad\, \times  \sum\limits_{m = 0}^\infty  \frac{1}{z^{m/\mu } }\sum\limits_{j = m }^{\left\lceil n/2 \right\rceil+2m} g_{\left\lceil n/2 \right\rceil+j,n} \left( \frac{r_0}{p_0^{1 /\mu } } \right)f_{\left\lceil n/2 \right\rceil+m,j - m}  .
\end{align*}
In passing to the last equality, we have used the fact that $g_{j,n}(x)= 0$ when $0 \le j \le n - 1$. The desired result follows by setting
\[
\widetilde f_{n,m} := \sum\limits_{j = m }^{\left\lceil n/2 \right\rceil+2m} g_{\left\lceil n/2 \right\rceil+j,n} \left( \frac{r_0}{p_0^{1 /\mu } } \right)f_{\left\lceil n/2 \right\rceil+m,j - m}.
\]

\section{Examples}\label{app}

In this section, we give two illustrative examples to demonstrate the applicability of our results.

\subsection{Confluent hypergeometric function $U(a,b,z)$ for large $b$ and $z$} The confluent hypergeometric functions are solutions of the second order linear ordinary differential equation
\[
zw''(z) + (b - z)w'(z) - aw(z) = 0.
\]
The standard solutions are denoted by $M(a,b,z)$ and $U(a,b,z)$ (see, for example, \cite[\S 13.2]{DLMF}). The $U$-function has the integral representation
\begin{equation}\label{eq43}
U(a,b,z) = \frac{1}{\Gamma (a)}\int_0^{ + \infty } e^{ - zs} s^{a - 1} (1 + s)^{b - a - 1} ds,
\end{equation}
valid for $\Re(a)>0$ and $\left|\arg z\right|<\frac{\pi}{2}$ \cite[eq. 13.4.4]{DLMF}. In this example, we are interested in the asymptotic behaviour of $U(a,b,z)$ when both $|b|$ and $|z|$ become large. Temme \cite[Sec. 10.4]{Temme2015} showed that if $a$ is fixed and $\Re(a)>0$, then
\begin{equation}\label{eq42}
U(a,b,\zeta b) = \frac{1}{\Gamma (a)}\zeta ^{1-b} (1 - \zeta )^{a - 1} e^{ - b(1-\zeta)} \sqrt {\frac{2\pi }{b}} \left( 1 + \mathcal{O}\left( \frac{1}{(1 - \zeta )^2 b} \right) \right),
\end{equation}
as $b\to+\infty$ provided $0 < \zeta  < 1$, and
\begin{equation}\label{eq40}
U(a,b,\zeta b) = \frac{1}{(\zeta  - 1)^a b^a }\left( 1 + \mathcal{O}\left( \frac{1}{(\zeta  - 1)^2 b} \right) \right),
\end{equation}
as $b\to+\infty$ provided $\zeta  > 1$. When $z \approx b$, we have
\begin{equation}\label{eq41}
U(a,b,b+\kappa) = \frac{\Gamma \left( \frac{a}{2} \right)}{2\Gamma (a)}\left( \frac{2}{b+\kappa} \right)^{a/2} \left( 1 + \mathcal{O}\left( \frac{\left|\kappa\right|+1}{b^{1/2} } \right) \right)
\end{equation}
as $b\to+\infty$, where $\kappa$ is an arbitrary fixed complex number. The proof of this result follows along the same lines as the proof of the special case $\kappa=-a-1$ given in the paper \cite{Lopez2010} by L\'opez and Pagola. The fractional powers in \eqref{eq42}--\eqref{eq41} assume their principal values.

The approximations \eqref{eq42} and \eqref{eq40} are suitable for numerical computation when $\left| z-b \right|$ is large compared with $b^{1/2}$, for otherwise the error terms of these approximations do not tend to zero as $b\to+\infty$. Similarly, the approximation \eqref{eq41} is useful only when $z =b+o\left(b^{1/2}\right)$. Thus, corresponding to values of $\left| z-b \right|$ that are comparable with $b^{1/2}$, there are gaps within which neither of these approximations is suitable. We call these gaps the transition regions.

\begin{figure}[!h]
\centering
\def\svgwidth{0.50\textwidth}
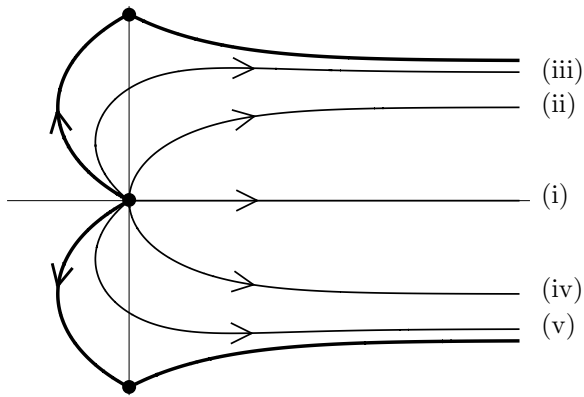
\caption{The steepest descent contour $\mathscr{P}\left(\theta\right)$ associated with the confluent hypergeometric function emanating from the origin when (i) $\theta=0$, (ii) $\theta=-\pi$, (iii) $\theta=-\frac{11\pi}{8}$, (iv) $\theta=\pi$ and (v) $\theta=\frac{11\pi}{8}$. The thick contours are the limiting paths $\mathscr{P}\left(\mp \frac{3\pi}{2}\right)$ through the points $\pm 2\pi i$.}
\label{figure1}
\end{figure}

In the forthcoming, we show how the main results of this paper can be used to derive an asymptotic expansion of the $U$-function which is valid in the transition regions. Let $\tau$ be any fixed complex number. Substituting $z=b+\tau b^{1/2}$ and $s=e^t-1$ in \eqref{eq43}, we obtain
\begin{equation}\label{eq48}
U(a,b,b+\tau b^{1/2}) = \frac{1}{\Gamma (a)}\int_0^{ + \infty } \exp \left( - bp(t) + b^{1/2} r(t) \right)q(t) dt,
\end{equation}
with
\[
p(t) = e^t  - t - 1,\quad q(t) = \left( e^t  - 1 \right)^{a - 1} e^{-at} ,\quad r(t) = \tau \left( 1 - e^t \right).
\]
The integral \eqref{eq48} is certainly convergent when $\left| b \right|^{1/2}  > \left| \tau  \right|\sec \theta$, $\Re(a)>0$ and $\left|\theta\right|<\frac{\pi}{2}$ with $\theta=\arg b$. The branch of $b^{1/2}$ has phase $\theta/2$. From the results of the paper \cite{Nemes2016}, we can assert that, as long as $\left|\theta\right|< \frac{3\pi}{2}$, the contour of integration in \eqref{eq48} can be deformed continuously into a steepest descent path $\mathscr{P}(\theta)$ of $e^{i\theta } p(t)$ without crossing any possible singularities of $q(t)$. Hence, by analytic continuation,
\begin{equation}\label{eq49}
U(a,b,b+\tau b^{1/2}) = \frac{1}{\Gamma (a)}\int_{\mathscr{P}(\theta)} \exp \left(  -b p(t) + b^{1/2} r(t) \right)q(t) dt 
\end{equation}
for sufficiently large values of $|b|$, provided that $\Re (a) > 0$ and $\left|\theta \right|<\frac{3\pi}{2}$. The contour $\mathscr{P}(\theta)$, for some specific values of $\theta$, is depicted in Figure \ref{figure1}. A further extension is possible by noting that the integral \eqref{eq49} remains convergent in the range $\frac{3\pi}{2} \leq \theta < 2\pi -\delta$ (respectively $-2\pi +\delta < \theta \leq -\frac{3\pi}{2}$) for any fixed positive $\delta$ if the integration is taken along the path $\mathscr{P}\left(\frac{3\pi}{2}-\frac{\delta}{2}\right)$ (respectively $\mathscr{P}\left(-\frac{3\pi}{2}+\frac{\delta}{2}\right)$) and $|b|$ is large enough. Hence, if we define $\mathscr{P}(\theta):=\mathscr{P}\left(\pm\frac{3\pi}{2}\mp\frac{\delta}{2}\right)$ for $\frac{3\pi}{2} \leq \pm\theta < 2\pi -\delta$, then \eqref{eq49} continues to hold in the wider range $\left|\theta\right|\leq 2\pi-\delta$ provided $|b|$ is sufficiently large.

We now apply Corollary \ref{thm2} to the integral \eqref{eq49} (with the extended definition of $\mathscr{P}(\theta)$), separately for each $\theta$ in the interval $-2\pi +\delta \leq \theta \leq 2\pi -\delta$. In the neighbourhood of $0$, the functions $p(t)$, $q(t)$ and $r(t)$ have convergent expansions
\[
p(t) = \frac{1}{2}t^2 + \frac{1}{6}t^3  + \cdots ,\quad q(t) = t^{a - 1} - \frac{a + 1}{2}t^a  + \cdots ,\quad r(t) = - \tau t - \frac{\tau}{2}t^2  -  \cdots .
\]
Thus, with the notation of Section \ref{intro}, $\mu=2$, $\lambda=a$ and $\nu=1$ (note that $\nu/\mu=1/2$ as required). A convenient choice for the domain $\textbf{T}$ is $\left\{ t:\left| \Im (t) \right| < 2\pi \right\}\setminus \left(  - \infty ,0 \right]$ (cf. Figure \ref{figure1}). It is readily checked that Conditions (i)--(vi) of Section \ref{intro} are all satisfied. Therefore, according to Corollary \ref{thm2},
\begin{gather}\label{eq50}
\begin{split}
U(a,b,b+\tau b^{1/2}) \sim \; & \frac{1}{\Gamma (a)}\frac{1}{b^{a/2}}\operatorname{Fi}\left( \frac{1}{2},\frac{a}{2}; - 2^{1/2} \tau \right)\sum\limits_{m = 0}^\infty  \frac{\widetilde f_{0,m}}{b^{m/2} }  \\ & + \frac{1}{\Gamma (a)}\frac{1}{b^{(a + 1)/2}}\operatorname{Fi}\left( \frac{1}{2},\frac{a + 1}{2}; - 2^{1/2} \tau \right)\sum\limits_{m = 0}^\infty  \frac{\widetilde f_{1,m} }{b^{m/2}} ,
\end{split}
\end{gather}
as $\left|b\right|\to+\infty$ uniformly in the sector $\left|\theta\right|\leq 2\pi-\delta$, with $a$ and $\tau$ being fixed and $\Re(a)>0$. The branch of $b^{a/2}$ is $\exp\left(a(\log|b|+i\theta)/2\right)$. Finally, the use of equation \eqref{eq52} leads us to the final form of the expansion in terms of the parabolic cylinder function:
\begin{gather}\label{eq51}
\begin{split}
U(a,b,b+\tau b^{1/2}) \sim \; & \frac{1}{b^{a/2}}\exp \left( \frac{\tau ^2 }{4} \right)U\left(a - \frac{1}{2},\tau \right)\sum\limits_{m = 0}^\infty  \frac{C_m (\tau,a )}{b^{m/2} }  \\ & + \frac{1}{b^{(a + 1)/2} }\exp \left( \frac{\tau ^2}{4} \right)U\left( a + \frac{1}{2},\tau \right)\sum\limits_{m = 0}^\infty  \frac{D_m (\tau,a )}{b^{m/2} } ,
\end{split}
\end{gather}
where the coefficients $C_m (\tau,a)$, $D_m (\tau,a )$ are polynomials in $\tau$ and $a$ given by
\[
C_m (\tau,a ) = 2^{1 - a/2} \widetilde f_{0,m} ,\quad D_m (\tau,a ) = 2^{1 - (a + 1)/2} a\widetilde f_{1,m} .
\]
Using formulae \eqref{eq20}, \eqref{eq21}, \eqref{eq22} and the procedure given in Section \ref{proof2}, we find the first few values of these coefficients to be
\begin{gather*}
C_0 (\tau,a) = 1,\quad C_1 (\tau,a ) =  - \frac{a}{3}\tau,\quad C_2 (\tau,a ) = \frac{a}{18}\tau ^4  - \frac{2a^2  + a}{12}\tau ^2  + \frac{8a^3  + 21a^2  + 10a}{36}, \\
C_3 (\tau,a ) =   - \frac{a}{162}\tau ^7  + \frac{2a^2  - 3a}{108}\tau ^5  - \frac{10a^3  - 45a^2  - 31a}{1620}\tau ^3  - \frac{40a^4  + 135a^3  + 188a^2  + 90a}{1620}\tau, \\
D_0 (\tau,a ) =  \frac{a}{3}\tau ^2  - \frac{2a^2  + 2a}{3},\quad D_1 (\tau,a ) =  - \frac{a}{18}\tau ^5  + \frac{4a^2  + a}{36}\tau ^3  - \frac{2a^2  + a}{12}\tau,\\
D_2 (\tau ,a) = \frac{a}{162}\tau ^8  - \frac{4a^2  - 11a}{324}\tau ^6  - \frac{5a^3  + 10a^2  + a}{270}\tau ^4  + \frac{100a^4  + 255a^3  + 227a^2  + 69a}{1620}\tau ^2  \\ \qquad \qquad \qquad\quad\;\;\,\, - \frac{40a^5  + 235a^4  + 443a^3  + 317a^2  + 69a}{810} .
\end{gather*}
We note that when $\left|\theta\right|\leq \frac{\pi}{2}-\delta$ and $\Re\left(\tau b^{1/2}\right)\geq 0$, the asymptotic expansion \eqref{eq51} is valid for arbitrary large values of $|\tau|$. This can be seen by applying Corollary \ref{thm2} to the integral representation \eqref{eq48} and taking into account remark (4) in Section \ref{intro}.

The expansion \eqref{eq51} was established under the condition that $\Re(a)>0$. We now show that this restriction can be removed. Suppose that $\Re(a)>-1$. We begin with the functional relation
\[
U(a,b,b + \tau b^{1/2} ) = (\tau b^{1/2}  + 2a+2)U(a+1,b,b + \tau b^{1/2} ) + (a+1)(b - a - 2)U(a + 2,b,b + \tau b^{1/2} )
\]
(cf. \cite[eq. 13.3.7]{DLMF}). Since $\Re(a)>-1$, we can employ the asymptotic expansion \eqref{eq51} to the right-hand side of this equality. In this way, we obtain an expansion for $U(a,b,b + \tau b^{1/2} )$ in terms of $U\left( a + \frac{1}{2},\tau \right)$, $U\left( a + \frac{3}{2},\tau \right)$ and $U\left( a + \frac{5}{2},\tau \right)$. The quantities $U\left( a + \frac{3}{2},\tau \right)$ and $U\left( a + \frac{5}{2},\tau \right)$ are expressible in terms of $U\left( a - \frac{1}{2},\tau \right)$ and $U\left( a + \frac{1}{2},\tau \right)$ by an application of the known recurrence relation for the parabolic cylinder function \cite[eq. 12.8.1]{DLMF}. An appeal to the uniqueness property of asymptotic series and analytic continuation then shows that \eqref{eq51} is valid for $\Re(a)>-1$. One then applies an induction argument to complete the proof.

As a closing remark, we mention that in a paper by Temme \cite{Temme1978}, an asymptotic expansion for $U(a,b,\zeta b)$ was established which is uniformly valid in any compact $\zeta$-interval of $(0,+\infty)$ (see also \cite[Subsec. 22.5.1]{Temme2015}). However, applying this asymptotic expansion in the neighbourhood of the point $z=b$ can be difficult because its coefficients possess a removable singularity at $\zeta=1$. In contrast, our expansion \eqref{eq51} has easily computable coefficients that do not involve a removable singularity at $z=b$ (where $\tau=0$).

\subsection{Associated Anger--Weber function $\textbf{A}_{-\rho}(z)$ for large $\rho$ and $z$} The functions of Anger and Weber are respectively defined by 
\[
\textbf{J}_{\rho}(z) := \frac{1}{\pi }\int_0^\pi \cos (\rho t - z\sin t)dt ,\quad \textbf{E}_\rho(z) := \frac{1}{\pi }\int_0^\pi \sin (\rho t - z\sin t)dt .
\]
Each is an entire function of $z$ and $\rho$. These functions appear as special solutions of inhomogeneous Bessel equations \cite[\S11.10(ii)]{DLMF}. The function $\textbf{J}_{\rho}(z)$ reduces to the Bessel function $J_n(z)$ when $\rho$ has the integer value $n$. When $\rho$ is not an integer, the two functions are distinct. Their difference can be expressed as $\textbf{J}_{\rho}(z)-J_\rho(z)=\sin(\pi\rho)\textbf{A}_\rho(z)$, where
\begin{equation}\label{eq57}
\textbf{A}_\rho  (z) := \frac{1}{\pi }\int_0^{ + \infty } e^{ - z\sinh t - \rho t} dt
\end{equation}
is the so-called associated Anger--Weber function. The integral \eqref{eq57} converges for any complex $\rho$ and $\left|\arg z\right|<\frac{\pi}{2}$ \cite[eq. 11.10.4]{DLMF}. In this example, we are concerned with the asymptotic behaviour of this function when $|\rho|$ and $|z|$ are both large. The function $\textbf{A}_\rho (z)$, with fixed $z/\rho>0$, has a simple asymptotic expansion as $\rho\to+\infty$ (see, e.g., \cite[eq. 11.11.8]{DLMF} or \cite[\S10.15]{Watson1944}). By contrast, the function $\textbf{A}_{-\rho } (z)$, for large positive values of $\rho$ and $z$, exhibits a more complicated behaviour. In this regard, Dingle \cite[p. 388]{Dingle1973} showed that
\begin{equation}\label{eq59}
\textbf{A}_{ - \rho } (\zeta \rho ) = \left( \frac{1 + \sqrt {1 - \zeta ^2 }}{\zeta } \right)^\rho  \frac{e^{ - \rho \sqrt {1 - \zeta ^2 } } }{\left( 1 - \zeta ^2 \right)^{1/4} }\sqrt {\frac{2}{\pi \rho}} \left( 1 + \mathcal{O}\left( \frac{1}{\left(1 - \zeta\right)^{3/2} \rho } \right)\right),
\end{equation}
as $\rho\to+\infty$ provided $0 < \zeta  < 1$. Watson \cite[\S10.15]{Watson1944} proved that
\begin{equation}\label{eq60}
\textbf{A}_{ - \rho } (\zeta \rho ) = \frac{1}{\pi (\zeta  - 1)\rho }\left( 1 + \mathcal{O}\left( \frac{1}{(\zeta  - 1)^3 \rho^2 } \right) \right),
\end{equation}
as $\rho\to+\infty$ with $\zeta  > 1$. For the case $z\approx \rho$, Airey \cite{Airey1918} gave
\begin{equation}\label{eq61}
\textbf{A}_{ - \rho } (\rho  + \kappa ) = \frac{2^{4/3}}{3^{7/6} \Gamma \left( \frac{2}{3} \right)(\rho  + \kappa )^{1/3} }\left( 1 + \mathcal{O}\left( \frac{\left| \kappa  \right| + 1}{\rho ^{1/3} } \right) \right)
\end{equation}
as $\rho\to+\infty$, where $\kappa$ is any fixed complex number (see also \cite[\S10.15]{Watson1944}). The fractional powers in \eqref{eq59}--\eqref{eq61} take their principal values.

The situation is similar to that encountered in the first example. The approximations \eqref{eq59} and \eqref{eq60} are of practical use only when $\left| z-\rho \right|$ is large compared with $\rho^{1/3}$, and the approximation \eqref{eq61} is suitable for numerical computation when $z =\rho+o\left(\rho^{1/3}\right)$. Therefore, there are transition regions between \eqref{eq59} and \eqref{eq61} and also between \eqref{eq60} and \eqref{eq61}, and in these transition regions $\left| z-\rho \right|$ is comparable with $\rho^{1/3}$. To bridge these gaps, Olver \cite[exer. 5.1, p. 336]{Olver1974} proposed the asymptotic formula
\begin{equation}\label{eq62}
\textbf{A}_{ - \rho } (\rho  + \tau \rho ^{1/3} ) = \frac{2^{1/3}}{\rho ^{1/3}}\operatorname{Hi}( - 2^{1/3} \tau ) + \mathcal{O}\left( \frac{1}{\rho} \right),
\end{equation}
valid for large positive values of $\rho$ and bounded real values of $\tau$. Here $\operatorname{Hi}(x)$ denotes Scorer's function (see Appendix \ref{appendixa}).

\begin{figure}[!ht]
\centering
\def\svgwidth{0.50\textwidth}
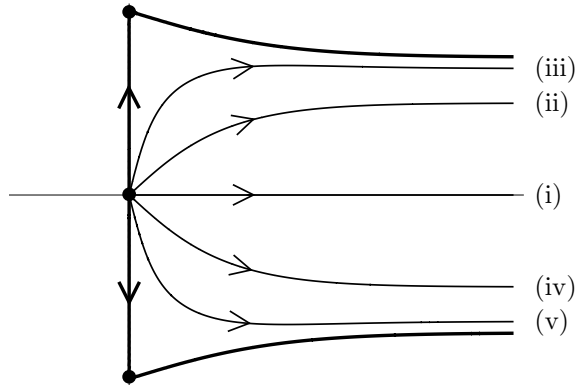
\caption{The steepest descent contour $\mathscr{P}\left(\theta\right)$ associated with the Anger--Weber function emanating from the origin when (i) $\theta=0$, (ii) $\theta=-\pi$, (iii) $\theta=-\frac{11\pi}{8}$, (iv) $\theta=\pi$ and (v) $\theta=\frac{11\pi}{8}$. The thick contours are the limiting paths $\mathscr{P}\left(\mp \frac{3\pi}{2}\right)$ through the points $\pm 2\pi i$.}
\label{figure2}
\end{figure}

In what follows, we apply the main results of this paper to extend \eqref{eq62} to a complete asymptotic expansion valid for complex values of $\rho$ and $\tau$ as $|\rho|\to+\infty$. We can proceed similarly to the case of the $U$-function. Let $\tau$ be an arbitrary fixed complex number. Substituting $z=\rho+\tau \rho^{1/3}$ in \eqref{eq57}, we obtain
\begin{equation}\label{eq56}
\textbf{A}_{ - \rho } (\rho  + \tau \rho ^{1/3} ) = \frac{1}{\pi }\int_0^{ + \infty } \exp \left(  - \rho p(t) + \rho ^{1/3} r(t) \right)q(t)dt ,
\end{equation}
with
\[
p(t) = \sinh t - t,\quad q(t) = 1,\quad r(t) =  - \tau \sinh t.
\]
The integral \eqref{eq56} is certainly convergent when $\left| \rho \right|^{2/3}  > \left| \tau  \right|\sec \theta$ and $\left|\theta\right|<\frac{\pi}{2}$ with $\theta=\arg \rho$. The branch of $\rho^{1/3}$ has phase $\theta/3$. According to the results of the paper \cite{Nemes2014}, the integration contour in \eqref{eq56} can be deformed continuously into a path of steepest descent $\mathscr{P}(\theta)$ of $e^{i\theta } p(t)$ when $\left|\theta\right|< \frac{3\pi}{2}$. Thus, by analytic continuation,
\begin{equation}\label{eq58}
\textbf{A}_{ - \rho } (\rho  + \tau \rho ^{1/3} ) = \frac{1}{\pi }\int_{\mathscr{P}(\theta)} \exp \left(  - \rho p(t) + \rho ^{1/3} r(t) \right)q(t)dt
\end{equation}
for all sufficiently large values of $|\rho|$, provided that $\left|\theta \right|<\frac{3\pi}{2}$. The contour $\mathscr{P}(\theta)$, for some specific values of $\theta$, is shown in Figure \ref{figure2}. Following the steps of the preceding subsection, we define $\mathscr{P}(\theta):=\mathscr{P}\left(\pm\frac{3\pi}{2}\mp\frac{\delta}{2}\right)$ for $\frac{3\pi}{2} \leq \pm\theta < 2\pi -\delta$, where $\delta$ is any fixed positive number. With this notation, \eqref{eq58} is valid for $\left|\theta\right|\leq 2\pi-\delta$ provided $|\rho|$ is large enough.

We now apply Theorem \ref{thm1} to the integral \eqref{eq58} (with the extended definition of $\mathscr{P}(\theta)$), separately for each $\theta$ in the interval $-2\pi +\delta \leq \theta \leq 2\pi -\delta$. The functions $p(t)$, $q(t)$ and $r(t)$ have power series expansions of the form
\[
p(t) = \frac{1}{6}t^3  + \frac{1}{120}t^5  +  \cdots ,\quad q(t) = 1,\quad r(t) =  - \tau t - \frac{\tau }{6}t^3  -  \cdots .
\]
Hence, with the notation of Section \ref{intro}, $\mu=3$, $\lambda=1$ and $\nu=1$ (note that $\nu/\mu=1/3$ as required). The domain $\textbf{T}$ can be taken to be the entire complex plane. It is easily checked that Conditions (i)--(vi) of Section \ref{intro} are all satisfied. Note that since $p(t)$ and $r(t)$ are odd functions and $q(t)$ is an even function, formulae \eqref{eq20} and \eqref{eq21} imply that $b_n=c_n=0$ for all even $n$. Consequently, $f_{n,m}  = 0$ when $n$ is odd. Thus, in this particular case, Theorem \ref{thm1} yields the asymptotic expansion
\begin{equation}\label{eq69}
\textbf{A}_{ - \rho } (\rho  + \tau \rho ^{1/3} ) \sim \frac{1}{\pi}\sum\limits_{n = 0}^\infty \frac{1}{\rho^{(2n + 1)/3} }\sum\limits_{m = 0}^{2n} f_{2n,m} \operatorname{Fi}\left( \frac{1}{3},\frac{2n + m + 1}{3}; - 6^{1/3} \tau  \right),
\end{equation}
as $\left|\rho\right|\to+\infty$ uniformly in the sector $\left|\theta\right|\leq 2\pi-\delta$, with $\tau$ being fixed. It is worth mentioning that a similar expansion was given by Dingle \cite[p. 250]{Dingle1973}. However, the analysis used in the derivation of his expansion is purely formal, and the region of validity is not examined. Using a method analogous to that of Section \ref{proof2}, we rearrange the asymptotic expansion \eqref{eq69} in the form
\[
\textbf{A}_{ - \rho } (\rho  + \tau \rho ^{1/3} ) \sim \frac{1}{\pi}\sum\limits_{n = 0}^2 \frac{1}{\rho ^{(2\left\lceil {n/2} \right\rceil  + 1)/3} }\operatorname{Fi}^{(n)} \left( \frac{1}{3},\frac{1}{3}; - 6^{1/3} \tau  \right)\sum\limits_{m = 0}^\infty \frac{\widehat f_{n,m}}{\rho ^{2m/3}},
\]
with
\begin{equation}\label{eq53}
\widehat f_{n,m}  = \sum\limits_{j = 2m}^{2\left\lceil {n/2} \right\rceil  + 4m} g_{2\left\lceil {n/2} \right\rceil  + j,n} ( - 6^{1/3} \tau )f_{2\left\lceil {n/2} \right\rceil  + 2m,j - 2m} .
\end{equation}
 Finally, with the aid of equations \eqref{eq54} and \eqref{eq65}, we obtain the desired form of the expansion in terms of the Scorer function and its derivative:
\begin{gather}\label{eq55}
\begin{split}
\textbf{A}_{ - \rho } (\rho  + \tau \rho ^{1/3} ) \sim \; & \frac{2^{1/3}}{\rho ^{1/3}}\operatorname{Hi}( - 2^{1/3} \tau )\sum\limits_{m = 0}^\infty  \frac{U_m (\tau )}{\rho ^{2m/3} }  \\ & + \frac{2^{2/3}}{\rho}\operatorname{Hi}'( - 2^{1/3} \tau )\sum\limits_{m = 0}^\infty  \frac{V_m (\tau )}{\rho ^{2m/3} }  + \frac{1}{\pi \rho }\sum\limits_{m = 0}^\infty  \frac{W_m(\tau) }{\rho ^{2m/3} }, 
\end{split}
\end{gather}
where the coefficients $U_m (\tau )$, $V_m (\tau )$ and $W_m (\tau )$ are polynomials in $\tau$ given by (with the convention that $\widehat f_{2,- 1}=0$)
\[
U_m (\tau ) = 2^{ - 1/3} 3^{2/3} \widehat f_{0,m} - \tau \widehat f_{2,m - 1} ,\quad V_m (\tau ) = 2^{ - 2/3} 3^{1/3} \widehat f_{1,m} ,\quad W_m (\tau ) = \widehat f_{2,m} .
\]
From \eqref{eq20}, \eqref{eq21}, \eqref{eq22}, \eqref{eq53} and the recurrence relations given in Section \ref{proof2}, the first few coefficients are found to be 
\begin{gather*}
U_0(\tau)=1,\quad U_1 (\tau ) =  - \frac{1}{5}\tau,\quad U_2 (\tau ) = -\frac{9}{100}\tau ^5+\frac{3}{35}\tau ^2,\quad U_3 (\tau ) = \frac{957}{7000}\tau ^6  - \frac{173}{3150}\tau ^3  - \frac{1}{225},\\ V_0 (\tau ) = \frac{3}{10}\tau ^2,\quad V_1 (\tau ) =  - \frac{17}{70}\tau ^3  + \frac{1}{70}, \quad V_2 (\tau ) =  - \frac{9}{1000}\tau ^7  + \frac{611}{3150}\tau ^4  - \frac{37}{3150}\tau, \\
W_0 (\tau ) = -\frac{1}{10}, \quad W_1 (\tau ) = \frac{9}{100}\tau^4+\frac{47}{700}\tau, \quad W_2(\tau)= - \frac{447}{3500}\tau ^5  -\frac{23}{600}\tau ^2.
\end{gather*}
The coefficients $U_m (\tau )$ and $V_m (\tau )$ seem to agree with those appearing in the analogous asymptotic expansion of the Bessel function $J_{\rho } (\rho  + \tau \rho ^{1/3} )$ \cite[\S10.19(iii)]{DLMF}. This is of no surprise considering that $\textbf{A}_{ - \rho } (z)$ is a solution of an inhomogeneous Bessel equation.

We remark that when $\left|\theta\right|\leq \frac{\pi}{2}-\delta$ and $\Re\left(\tau \rho^{1/3}\right)\geq 0$, the asymptotic expansion \eqref{eq55} is valid for arbitrary large values of $|\tau|$. This can be shown by employing Theorem \ref{thm1} to the integral representation \eqref{eq56} and taking into account remark (4) in Section \ref{intro}.

We close this subsection by mentioning that asymptotic expansions of the function $\textbf{A}_{ - \rho } (\zeta\rho)$, for large $\rho$, which are uniformly valid for $\zeta$ in the intervals $0<\zeta\leq 1$ and $\zeta\geq 1$, respectively, can be found in the book of Olver \cite[pp. 352--360]{Olver1974}. However, in contrast to our expansion \eqref{eq55}, the evaluation of the coefficients of these asymptotic expansions near the point $z=\rho$ is a difficult problem.

\section{Conclusion}\label{disc}

In this paper, we have stated and proved theorems concerning the asymptotic expansion of contour integrals of the form
\begin{equation}\label{eq67}
\int_a^b  \exp \left(  - zp(t) + z^{\nu /\mu } r(t) \right)q(t)dt,
\end{equation}
in which $z$ is a large real or complex parameter, $p(t)$, $q(t)$ and $r(t)$ are analytic functions of $t$, and the positive constants $\mu$ and $\nu$ are related to the local behaviour of the functions $p(t)$ and $r(t)$ near the endpoint $a$.

The main theorem of the paper includes as special cases several important asymptotic methods for integrals. In particular, if $r(t)$ is identically zero, our main result coincides with that of Olver \cite{Olver1970}. The paper \cite{Olver1970} also supplies explicit bounds for the error terms associated with the asymptotic expansions. It may be possible to obtain similar error bounds in the current more general context, but this has not been explored in the present paper.

The new results have potential applications in the asymptotic theory of special functions in transition regions, and we have illustrated this by two examples.

Asymptotic expansions similar to ours for integrals of type \eqref{eq67} and also for solutions of homogeneous and inhomogeneous linear ordinary differential equations are given in the book by Dingle \cite{Dingle1973}. Dingle's derivation of his expansions is based on non-rigorous methods. The analysis that we have given in this paper provides a rigorous justification for some of his results, especially those discussed in Chapters X and XI (uniform asymptotic expansions for integrals). Most results concerning differential equations (Chapters XV, XVII and XX), however, still lack rigorous mathematical foundations and their investigation should be the subject of future research.

\appendix

\section{Basic properties of Fax\'{e}n's integral}\label{appendixa}

In this appendix, we discuss some of the basic properties of the Fax\'{e}n integral defined by equation \eqref{eq35}. It is assumed throughout that $0 \le \alpha  < 1$ and $\Re (\beta ) > 0$.

By uniform convergence we may expand the factor $\exp \left( xt^\alpha \right)$ in ascending powers of $x t^\alpha$ and integrate \eqref{eq35} term by term. In this way we obtain the following expansion, valid for all $x$: 
\[
\operatorname{Fi}(\alpha ,\beta ;x) = \sum\limits_{n = 0}^\infty  \Gamma (\alpha n + \beta )\frac{x^n}{n!} .
\]
This expansion is useful for computing $\operatorname{Fi}(\alpha ,\beta ;x)$ when $x$ is small or moderate in size. For large values of $|x|$ asymptotic expansions should be used instead. Fax\'{e}n's integral is a particular case of the so-called Wright function since $\operatorname{Fi}(\alpha ,\beta ;x) = {}_1\Psi _1 \left( {\alpha ,\beta ;0,1;x} \right)$. Accordingly, the large-$x$ asymptotic behaviour of $\operatorname{Fi}(\alpha ,\beta ;x)$ can be determined by employing the well-established asymptotic theory of the Wright function \cite{Paris2014}. In particular, to leading order,
\[
\operatorname{Fi}(\alpha ,\beta ;x) \sim \left( {\alpha x} \right)^{(2\beta  - 1)/(2 - 2\alpha )} \exp \left( {\left( {1 - \alpha } \right)\left( {\alpha ^\alpha  x} \right)^{1/(1 - \alpha )} } \right)\sqrt {\frac{{2\pi }}{{1 - \alpha }}} 
\]
and
\[
\operatorname{Fi}(\alpha ,\beta ;-x) \sim \Gamma \left( {\frac{\beta }{\alpha }} \right)\frac{1}{{\alpha x^{\beta /\alpha } }}
\]
as $x\to +\infty$. For further representations, we refer the reader to \cite{Kaminski1997}.

In the remaining part of this appendix, we show that for certain specific values of the parameters $\alpha$ and $\beta$, the function $\operatorname{Fi}(\alpha ,\beta ;x)$ is expressible in terms of well-known special functions. These representations lead to simplified forms of the asymptotic expansions derived in Section \ref{app}. When $\alpha=\frac{1}{2}$ and $\beta$ is any complex number with positive real part, Fax\'{e}n's integral is essentially a parabolic cylinder function \cite[\S12.2]{DLMF}. Indeed, substituting $t=\frac{1}{2} s^2$, we obtain
\begin{gather}\label{eq52}
\begin{split}
\operatorname{Fi}\left( \frac{1}{2},\beta ;x \right) & = \int_0^{ + \infty } \exp \left(  - t + xt^{1/2} \right)t^{\beta  - 1} dt  = 2^{1 - \beta } \int_0^{ + \infty } \exp \left(  - \frac{1}{2}s^2 + \frac{x}{2^{1/2}}s \right)s^{2\beta  - 1} ds  \\ & = 2^{1 - \beta } \Gamma (2\beta )\exp \left( \frac{x^2}{8} \right)U\left( 2\beta  - \frac{1}{2}, - 2^{ - 1/2} x \right) ;
\end{split}
\end{gather}
compare \cite[eq. 12.5.1]{DLMF}. When $\alpha=\beta= \frac{1}{3}$, the substitution $t = \frac{1}{3}s^3$ produces
\begin{gather}\label{eq54}
\begin{split}
\operatorname{Fi}\left( \frac{1}{3},\frac{1}{3};x \right) & = \int_0^{ + \infty } \exp \left(  - t + xt^{1/3}  \right)t^{ - 2/3} dt \\ & = 3^{2/3} \int_0^{ + \infty } \exp \left(  - \frac{1}{3}s^3  + \frac{x}{3^{1/3}}s \right)ds  = 3^{2/3} \pi \operatorname{Hi}( 3^{ - 1/3} x )
\end{split}
\end{gather}
(cf. \cite[eq. 9.12.20]{DLMF}), where $\operatorname{Hi}(x)$ is Scorer's function \cite[\S9.12]{DLMF}. The function $\operatorname{Hi}(x)$ is known to satisfy the inhomogeneous Airy equation
\begin{equation}\label{eq65}
\operatorname{Hi}''(x) - x\operatorname{Hi}(x) = \frac{1}{\pi}.
\end{equation}
For other special cases, see \cite[p. 332]{Olver1974}.

\section{}\label{appendixb}

In this appendix, we prove the estimate \eqref{eq8}. With the notation of Section \ref{proof}, we define
\[
F(w,y) := v^{1 - \lambda /\mu } f(v,y/v)
\]
for sufficiently small $|w|$ and $\left|\arg y\right|<\frac{\pi}{2}$. Thus, by \eqref{eq29}, $F(w,y)$ is an analytic function of $w$ and it has the power series expansion
\[
F(w,y) = \sum\limits_{n = 0}^{\infty} \left( \sum\limits_{m = 0}^n f_{n,m} y^{m\nu /\mu }  \right)w^n,
\]
valid within some circle $\left|w\right|\leq\rho$, where $\rho$ is independent of $y$. The branch of $y^{\nu/\mu}$ has phase $\nu \arg y/\mu$. We truncate this power series after $N\geq 0$ terms and write
\[
F(w,y) = \sum\limits_{n = 0}^{N - 1} \left( \sum\limits_{m = 0}^n f_{n,m} y^{m\nu /\mu }  \right)w^n  + w^N F_N (w,y),
\]
with the remainder
\begin{equation}\label{eq15}
F_N (w,y) = f_N (v,y/v)
\end{equation}
(compare \eqref{eq30}). By choosing $k$ sufficiently close to $a$, we can assume that $\left|p(k)-p(a)\right|^{1/\mu}<\rho$. Then
\begin{equation}\label{eq16}
\left| F_N (w,y) \right| \le \frac{\rho ^{ - N} }{2\pi }\frac{M(\rho )}{\rho  - \left| w \right|},\qquad M(\rho ) = \mathop {\max }\limits_{\left| \zeta  \right| = \rho } \left| F(\zeta ,y) \right|,
\end{equation}
for any $w$ satisfying $w \leq \left|p(k)-p(a)\right|^{1/\mu}$ (see, for example, \cite[p. 70]{Markushevich1965}). It is readily seen from \eqref{eq14} that
\begin{equation}\label{eq17}
F(\zeta ,y) = \mathcal{O}(1)\exp \left( C_k y^{\nu /\mu } \right),\quad \left|\zeta\right|=\rho,
\end{equation}
where the implied constant depends only on $k$ and $\rho$, and $C_k$ is an assignable constant which may depend on $k$. Hence, from \eqref{eq5}, \eqref{eq15}, \eqref{eq16} and \eqref{eq17}, we can infer that
\[
f_N (v,z) = \mathcal{O}(1)\exp \left( C_k (zv)^{\nu /\mu } \right),
\]
uniformly with respect to $z$ and $v$, if $v$ lies on the segment connecting $0$ to $\kappa$.

\section{}\label{appendixc}

In this appendix, we show that if $0\leq \alpha<1$, and $\beta$ and $x$ are arbitrary (fixed) complex numbers, then
\begin{equation}\label{eq27}
\int_\zeta ^\infty \exp \left(  - t + xt^\alpha \right)t^{\beta  - 1} dt = \mathcal{O}(1)\exp \left(  - \zeta  + x\zeta ^\alpha \right)\zeta ^{\beta  - 1},
\end{equation}
uniformly in the region $\left|\arg \zeta\right| \leq \frac{\pi}{2}$ as $\left|\zeta\right|\to+\infty$. To this end, we perform a change of integration  variable from $t$ to $u$ by $t = \zeta (u + 1)$ in \eqref{eq27}, and deform the contour of integration so that $\arg u = -\arg \zeta$ along the new contour. Thus, the left-hand side of \eqref{eq27} becomes
\begin{gather}\label{eq28}
\begin{split}
& \zeta ^\beta  \int_0^{\infty e^{ - i\arg \zeta } } \exp \left(  - \zeta \left( u + 1\right) + x\left( \zeta \left( u + 1 \right) \right)^\alpha \right)\left( u + 1 \right)^{\beta  - 1} du \\ & = \exp \left(  - \zeta  + x\zeta ^\alpha  \right)\zeta ^\beta  \int_0^{\infty e^{ - i\arg \zeta } } \exp \left(  - \zeta u + x\zeta ^\alpha  \left( \left( u + 1 \right)^\alpha   - 1\right) \right)\left( u + 1 \right)^{\beta  - 1} du ,
\end{split}
\end{gather}
for any $\zeta$ in the closed sector $\left|\arg \zeta\right| \leq \frac{\pi}{2}$. Since $0\leq \alpha<1$, we may assert that
\[
\Re \left(  - \zeta u + x\zeta ^\alpha  \left( \left( u + 1 \right)^\alpha   - 1 \right) \right) =  - \left| \zeta  \right|\left| u \right|\Re \left( 1 - \frac{x}{\zeta ^{1 - \alpha } }\frac{\left( u + 1 \right)^\alpha   - 1}{u} \right) \le  - \frac{1}{2}\left| \zeta  \right|\left| u \right|,
\]
for sufficiently large $\left|\zeta\right|$. Therefore, for large $\left|\zeta\right|$, the integral in the second line of \eqref{eq28} can be estimated as follows:
\begin{multline*}
\left| \int_0^{\infty e^{ - i\arg \zeta } } \exp \left(  - \zeta u + x\zeta ^\alpha  \left( \left( u + 1 \right)^\alpha   - 1 \right) \right)\left( u + 1 \right)^{\beta  - 1} du \right|  \\ \le \int_0^{\infty e^{ - i\arg \zeta } } e^{ - \frac{1}{2}\left| \zeta  \right|\left| u \right|} \left| \left( u + 1 \right)^{\beta  - 1}  \right| \left| du \right| = \int_0^{ + \infty } e^{ - \frac{1}{2}\left| \zeta  \right|s} \left| \left( se^{ - i\arg \zeta }  + 1 \right)^{\beta  - 1} \right|ds,
\end{multline*}
with $s=\left| u \right|$. Since the factor $\left( se^{ - i\arg \zeta }  + 1 \right)^{\beta  - 1}$ grows sub-exponentially in $s$, the last quantity is $\mathcal{O}(\left|\zeta\right|^{-1})$, which completes the proof of the estimate \eqref{eq27}.

\end{document}